\documentclass{aa}
\usepackage{graphicx}
\usepackage{txfonts}
\usepackage[utf8]{inputenc}
\usepackage{bm}
\usepackage{tikz}
\usepackage{algorithm}
\usepackage{algorithmicx}
\usepackage{algpseudocode} 
\usepackage{listings}
\usepackage{longtable}
\usepackage{booktabs}
\usepackage{tabularx}
\usepackage{svg}
\usepackage[colorlinks=true, allcolors=blue]{hyperref}
\usepackage{natbib}
\bibpunct{(}{)}{;}{a}{}{,} 
\usepackage{xcolor}
\usepackage{siunitx}
\begin{document}

\title{An adaptive symplectic integrator \\ for gravitational dynamics}

\author{Keqi Ye\inst{1,2} \and Zizhe Cai\inst{1,2} \and Mingji Wang\inst{3} \and Kun Yang\inst{1,2} \and Xiaodong Liu\inst{1,2}\thanks{Corresponding author: liuxd36@mail.sysu.edu.cn}}

\institute{
    School of Aeronautics and Astronautics, Shenzhen Campus of Sun Yat-sen University, Shenzhen, Guangdong 518107, China
\and
    Shenzhen Key Laboratory of Intelligent Microsatellite Constellation, Shenzhen Campus of Sun Yat-sen University, Shenzhen, Guangdong 518107, China
\and
    School of Mechanics and Aerospace Engineering, State Key Laboratory of Structural Analysis, Optimization and CAE Software for Industrial Equipment, Dalian University of Technology, Dalian, Liaoning 116081, China
}

\abstract
{This paper presents an adaptive symplectic integrator, SQQ-PTQ, developed on the basis of the fixed-step symplectic integrator SQQ. To mitigate the Runge phenomenon, SQQ-PTQ employs Chebyshev interpolation for approximating the action, enhancing both the precision and stability of the interpolation. In addition, to reduce the computational cost of evaluating interpolation functions, SQQ-PTQ introduces a projection method that improves the efficiency of these computations. A key feature of SQQ-PTQ is its use of the time transformation to implement an adaptive time step. To address the challenge of computing complicated Jacobian matrices attributed to the time transformation, SQQ-PTQ adopts a quasi-Newton method based on Broyden’s method. This strategy accelerates the solution of nonlinear equations, thereby improving the overall computational performance. The effectiveness and robustness of SQQ-PTQ are demonstrated via three numerical experiments. In particular, SQQ-PTQ demonstrates adaptability in handling close-encounter problems. Moreover, during long-term integrations, SQQ-PTQ maintains the energy conservation, further confirming its advantages as a symplectic algorithm.}


   \maketitle

\section{Introduction}
\label{intro}
Symplecticity is characteristic of Hamiltonian systems and numerical results obtained by symplectic integrators maintain long-term stability compared to nonsymplectic integrators.  Thus, symplectic integrators have proven to be powerful tools in studies of gravitational dynamics \citep{prefer1,prefer2,prefer3,prefer4}.

An important approach to constructing symplectic integrators is based on the use of generating functions \citep{gf1,gf2,gf3,gf4}. \citet{SQQ} constructed a high-order and fixed-step symplectic integrator SQQ for Hamiltonian systems based on the principle of least action and generating function. In the acronym SQQ, S denotes the symplectic nature of  the integrator, while the two Qs indicate that generalized coordinates are used as independent variables at both sides in the generating function \citep{SQQ}. Overall, SQQ is a general-purpose integrator, and has been applied to trajectory planning \citep{peng2023symplectic}, kinodynamic planning, optimal control of cable-driven tensegrity manipulators \citep{song2020novel}, and multibody dynamics problems \citep{phj1,phj2,lifei1}. The numerical results demonstrate that SQQ exhibits high accuracy and effective energy conservation behavior.

However, in the numerical computation of multi-body gravitational problems, traditional fixed-step symplectic algorithms face two primary challenges. The first challenge arises from the computational efficiency issues caused by the multi-scale nature of time in such systems. For systems with strong variations in interaction, such as high-eccentricity or hyperbolic Kepler problems, the fixed-step symplectic integrator becomes inefficient. In contrast, adaptive time-step integrators can handle these scenarios effectively by employing smaller time steps in rapidly changing regions and larger time steps in slowly varying regions \citep{dehnen2011n,springel2005cosmological}. Unfortunately, traditional step-size control strategies destroy the long-term conservation properties of symplectic integrators \citep{preto1999class,lose_sym2,blanes2012explicit}.

One approach to constructing adaptive time-step symplectic integrators involves the application of time transformation, which decouples the time step from the integration step \citep{hairer1997,mikkola1999algorithmic,preto1999class,time_tran_1,sym_timetran2}. However, time transformation often leads to a non-separable Hamiltonian system, which, in turn, transforms an explicit algorithm into an implicit one \citep{hairer1997,huang1997adaptive}. In general, implicit algorithms require solutions to nonlinear equations. Due to the complexity introduced by the time-transformed Hamiltonian, computing the Jacobian matrix for solving this nonlinear system can become computationally expensive. This leads to the second challenge: how to efficiently solve the nonlinear equations at each time step. To address this challenge, we employed a quasi-Newton method based on Broyden's method \citep{broyden1965class} to accelerate the inversion of the Jacobian matrix, which has been widely used in many research fields \citep[e.g.][]{bulin2021efficient,chen2022efficient,sherman1950adjustment}. Broyden’s method iteratively updates an approximation to the Jacobian inverse during the process of solving the nonlinear equations. By avoiding the explicit computation of the Jacobian matrix, this approach significantly reduces computational costs and improves efficiency.

To address the two important challenges mentioned above, our proposed method combines the time transformation and the quasi-Newton approach, offering a robust framework for efficient and accurate simulations in multi-body gravitational problems. This paper is organized as follows. Section \ref{sec:2new} provides the fundamental theory of SQQ based on Chebyshev interpolation method. Section \ref{sec:NumImpl} presents the numerical implementation, which includes three parts. The first part introduces a projection method that accelerates the computation of interpolation functions and their derivatives without requiring matrix inversion. The second part introduces the concept of the time transformation and discusses the selection of the step-size control function. The last part presents an efficient approach for solving nonlinear equations. Finally, Section \ref{sec:ne} includes several numerical examples to demonstrate the effectiveness of the proposed method.

\section{The symplectic algorithm based on Generating Function}
\label{sec:2new}
This section gives a brief introduction to SQQ (for more details, readers are referred to \citet{SQQ}). The action of a conservative system is the integral of the Lagrangian function over time
\begin{equation}
    S = \int_{t_a}^{t_b} \left[ \mathbf{p}^T \dot{\mathbf{q}} - H(\mathbf{q}, \mathbf{p}) \right] dt
    \label{S}
,\end{equation}
where $\mathbf{q}$ is the generalized displacement, \(\dot{\mathbf{q}}\) is the generalized velocity, and \(\mathbf{p}\) is the generalized momentum. The principle of least action is given by
\begin{equation}
\begin{aligned}
\delta S = & \int_{t_a}^{t_b} \left[ (\delta \dot{\mathbf{q}})^T \mathbf{p} - (\delta \mathbf{q})^T \frac{\partial H}{\partial \mathbf{q}} \right] dt + \int_{t_a}^{t_b} (\delta \mathbf{p})^T \left( \dot{\mathbf{q}} - \frac{\partial H}{\partial \mathbf{p}} \right) dt \\
= & \int_{t_a}^{t_b} (\delta \mathbf{q})^T \left[ -\dot{\mathbf{p}} - \frac{\partial H}{\partial \mathbf{q}} \right] dt  + \int_{t_a}^{t_b} (\delta \mathbf{p})^T \left( \dot{\mathbf{q}} - \frac{\partial H}{\partial \mathbf{p}} \right) dt \\
& + \mathbf{p}^T \delta \mathbf{q}\bigg|_{t_a}^{t_b} = 0.
\label{Eq:1}
\end{aligned}
\end{equation}
Equation~\ref{Eq:1} shows that if the Hamiltonian canonical differential equation is satisfied within the time step, the action will satisfy the following equation,
\begin{equation}
\delta S = \mathbf{p}^T(t_b) \, \mathrm{d}\mathbf{q}(t_b) - \mathbf{p}^T(t_a) \, \mathrm{d}\mathbf{q}(t_a).
\label{Eq:2}
\end{equation}
The physical meaning of Eq.~\ref{Eq:2} is that the action, \(S\), depends only on the generalized displacement at both ends of the time interval \(\left[ {{t_a},{t_b}} \right]\). In particular, \(S\) is a function of \(\mathbf{q}\left( {{t_a}} \right)\) and \(\mathbf{q}\left( {{t_b}} \right)\). From Eq.~\ref{Eq:2}, we get
\begin{equation}
    \mathbf{p}\left( t_a \right) =  - \frac{\partial S\left( \mathbf{q}\left( t_a \right), \mathbf{q}\left( t_b \right) \right)}{\partial \mathbf{q}\left( t_a \right)},\quad 
    \mathbf{p}\left( t_b \right) = \frac{\partial S\left( \mathbf{q}\left( t_a \right), \mathbf{q}\left( t_b \right) \right)}{\partial \mathbf{q}\left( t_b \right)}
\label{Eq:6}
.\end{equation}
Once the generating function is determined, the generalized displacement \(\mathbf{q}\left( t_b \right)\) and the generalized momentum \(\mathbf{p}\left( t_b \right)\) can be calculated according to Eq.~\ref{Eq:6}. In this context, SQQ employs the Lagrange interpolation with equidistant nodes and then the Gaussian quadrature is used to calculate the approximate value of \( S \). However, a higher order interpolation with equidistant nodes often leads to the Runge phenomenon. To avoid this phenomenon, we used the Chebyshev interpolation in this work \citep{cai2012chebyshev,lanczos1938trigonometric}.

The generalized coordinate and momentum vectors as functions of time, \(t\), are expressed as:
\begin{equation}
\mathbf{q}(t) = \sum_{k=0}^{m} M_k(t) \mathbf{q}_k, \quad \mathbf{p}(t) = \sum_{k=0}^{n} N_k(t) \mathbf{p}_k,
\label{Eq:7}
\end{equation}
where \(M_k(t)\) and \(N_k(t)\) are the interpolation basis functions, while \(\mathbf{q}_k\) and \(\mathbf{p}_k\) are defined as:
\begin{equation}
\mathbf{q}_k = \left\{ q_k^1, q_k^2, \cdots, q_k^d \right\}^\mathrm{T}, \quad k = 0, 1, \cdots, m
,\end{equation}
\begin{equation}
\mathbf{p}_k = \left\{ p_k^1, p_k^2, \cdots, p_k^d \right\}^\mathrm{T}, \quad k = 0, 1, \cdots, n
,\end{equation}
with \(d\) representing the degrees of freedom (DOF) of the system. Here, \(m\) and \(n\) denote the number of interpolation points for the generalized displacement and generalized momentum, respectively.

By substituting Eq.~\ref{Eq:7} into Eq.~\ref{S} and then combining it with Eq.~\ref{Eq:6}, we obtain
\begin{equation}
\begin{aligned}
&\frac{\partial S}{\partial \mathbf{q}_0} + \mathbf{p}_0  = 0, \\
&\frac{\partial S}{\partial \mathbf{q}_i} = 0, \quad i = 1, 2, \cdots, m-1, \\
&\frac{\partial S}{\partial \mathbf{p}_i} = 0, \quad i = 0, 1, \cdots, n.
\end{aligned}
\label{Eq:14a}
\end{equation}
and
\begin{equation}
\mathbf{p}_n = \frac{\partial S}{\partial \mathbf{q}_m}
\label{Eq:14b}
.\end{equation}

Generally, Eqs.~\ref{Eq:14a} is composed of nonlinear algebraic equations. Due to the introduction of the step-size control function in this paper (see Section \ref{Time Transformation}), the analytical derivation and calculation of the Jacobian matrix become complicated and challenging. This paper employs a quasi-Newton method and accelerates the iterative computation using Broyden's method (see Section \ref{qnbaesdbm}).

\section{Numerical implementation}\label{sec:NumImpl}
\subsection{Efficient computation of interpolation functions using the projection method}
To solve Eq.~\ref{Eq:14a}, the action \( S \) was approximated using Gaussian integration, which requires the calculation of Gaussian points and their corresponding interpolation functions, \( M_k(t) \) and \( N_k(t) \), at each time step. Since this process involves matrix inversion, it is computationally expensive. However, because the number of interpolation points remains constant during the computation, these functions do not need to be recomputed at each step. Instead, they can be efficiently determined using a projection method. This section will illustrate this method by using \( N_k(t) \) and \( \dot{N}_k(t) \) as examples.

Thus, we let Gaussian integral points and their weights in the interval $\left[ {-1, 1} \right]$ be
\begin{equation}
{\hat \xi _j}, \, {\hat \omega _j},\quad j = 1, 2, \cdots, g
,\end{equation}
where $g$ represents the number of Gaussian quadrature points, and ${\bm{\hat {N}}}$ and ${\bm{\hat {\dot{ N}}}}$ are the interpolation function and its derivative, respectively. In any time interval $\left[ {{t_a}, {t_b}} \right]$ ($0 \le {t_a} < {t_b}$), the Gaussian integral points and their weights can be obtained as follows:
\begin{equation}
{\xi _j} = {t_a} + \frac{{{\hat \xi _j} + 1}}{2}L,\quad {\omega _j} = \frac{L}{2}{\hat \omega _j},\quad j = 1, 2, \cdots, g
\label{Eq:17}
,\end{equation}
where \( L = t_b - t_a \). As illustrated in Fig.~\ref{fig:intervals}, \(\hat{\xi}\) and \(\xi\) (\(\xi \neq \hat{\xi}\)) denote the Gaussian quadrature points on the intervals \(\left[ -1, 1 \right]\) and \(\left[ t_a, t_b \right]\), respectively. These points share equivalent "length coordinates" denoted by \(\lambda_1\) and \(\lambda_2\), which represent the proportion of the interval's total length occupied by a given coordinate point. Consequently, the interpolation function within the interval \(\left[ t_a, t_b \right]\) can be expressed as
\begin{equation}
{\bm{N}}\left( {{{ \xi }_j}} \right) = {\bm{\hat N}}\left( {\hat \xi}_j \right),\quad j = 1, 2, \cdots, g
\label{Eq:18}
,\end{equation}
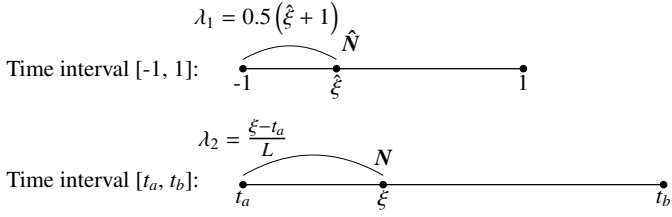
\begin{figure}
\centering
\resizebox{\columnwidth}{!}{ 
\begin{tikzpicture}[every node/.append style={font=\fontsize{14pt}{16pt}\selectfont}]
    \draw[thick] (0,0) -- (6,0);
    \filldraw (0,0) circle (2pt) node[below] {-1};
    \filldraw (6,0) circle (2pt) node[below] {1};
    \filldraw (2,0) circle (2pt) node[below] {$\hat{\xi}$};
    \node[above] at (2.3, 0.25) {$\bm{\hat{N}}$}; 
    \node[above] at (-3,-0.35) {Time interval [-1, 1]:}; 
    \node[above] at (0.5,0.6) {${\lambda _1} = 0.5\left( {\hat \xi + 1} \right)$};
    \draw (0,0.2) to [bend left=30] (2,0.2);
    \draw[thick] (0,-2.5) -- (9,-2.5);
    \filldraw (0,-2.5) circle (2pt) node[below] {$t_a$};
    \filldraw (9,-2.5) circle (2pt) node[below] {$t_b$};
    \filldraw (3, -2.5) circle (2pt) node[below] {$\xi$};
    \node[above] at (3, -2.2) {$\bm{N}$}; 
    \node[above] at (-3,-2.75) {Time interval [$t_a$, $t_b$]:};
    \node[above] at (0.0,-2.0) {${\lambda _2} = \scalebox{1.25}{$\frac{{\xi - {t_a}}}{L}$}$};
    \draw (0,-2.3) to [bend left=30] (3, -2.3);
\end{tikzpicture}
}
\caption{Illustration of the length coordinates for the time intervals $\left[ {-1, 1} \right]$ and $\left[ {{t_a}, {t_b}} \right]$. Although the horizontal coordinates of $\hat{\xi}$ and $\xi$ are not equal and their length coordinates $\lambda _1$ and $\lambda _2$ are equal, the values of the corresponding interpolation functions are equal (Eq.~\ref{Eq:18}).}
\label{fig:intervals}
\end{figure}
Taking the derivative of both sides of the Eq.~\ref{Eq:18}, we obtain
\begin{equation}
{\bm{\dot N}}\left( {{\xi _j}} \right) = \frac{2}{L}{\bm{\hat{\dot N}}}\left( {\hat \xi_j} \right),\quad j = 1, 2, \cdots, g
\label{Eq:19}
.\end{equation}

Clearly, it is only necessary to compute ${\bm{\hat N}}$ and ${\bm{\hat {\dot N}}}$, and then use the Eqs.~\ref{Eq:18} and \ref{Eq:19} to obtain ${\bm{N}}$ and ${\bm{\dot N}}$ in any time interval. The variables ${\bm{\hat N}}$ and ${\bm{\hat {\dot N}}}$ can be calculated using the following equations
\begin{equation}
{\bm{\hat N}}\left( {t} \right) = {\bm{G}}^{ - 1}{\bm{p}}\left( {t} \right), \quad {\bm{\dot N}}\left( {t} \right) = {\bm{G}}^{ - 1}{\bm{\dot p}}\left( {t} \right)
\label{Eq:14}
,\end{equation}
where
\begin{equation}
{\bm{G}} = \left[ {\begin{array}{*{20}{l}}
1&{{t_1}}&{t_1^2}& \cdots &{t_1^{n - 1}}\\
1&{{t_2}}&{t_2^2}& \cdots &{t_2^{n - 1}}\\
 \vdots & \cdots & \cdots & \cdots & \vdots \\
 \vdots & \cdots & \cdots & \cdots & \vdots \\
1&{{t_n}}&{t_n^2}& \cdots &{t_n^{n - 1}}
\end{array}} \right],\, {\bm{p}} = \left[ {\begin{array}{*{20}{l}}
1\\
{{t}}\\
{t^2}\\
 \vdots \\
{t^n}
\end{array}} \right],{\bm{\dot p}} = \left[ {\begin{array}{*{20}{l}}
0\\
1\\
{2{t}}\\
 \vdots \\
{nt^{n - 1}}
\end{array}} \right]
\end{equation}
and ${t_1}, {t_2}, \cdots, {t_n}$ are the Chebyshev interpolation points on the interval $\left[ {-1, 1} \right]$.

In this paper, a warm-up procedure is used to compute ${\bm{\hat N}}( {\hat \xi})$ and ${\bm{\hat{\dot N}}}( {\hat \xi} )$ in advance. Vectors ${\bm{N}}\left( {\xi } \right)$ and ${\bm{\dot N}}\left( {\xi } \right)$ are calculated at each step by using Eqs.~\ref{Eq:18} and \ref{Eq:19}, which reduces the computational costs of the matrix inversion. Similarly, this approach is also applied to compute ${\bm{M}}( { \xi})$ and ${\bm{{\dot M}}}( {\xi} )$.
\subsection{Time transformation}
\label{Time Transformation}
The time transformation \(t \leftrightarrow \tau\) can be employed to implement variable step sizes. Integrating the original system with a variable step size \(\Delta t\) is equivalent to integrating the transformed system with a constant step size, \(\Delta \tau\). The relationship between the transformed time, \(\tau\), and the original time, \(t\), is given by
\begin{equation}
    \frac{\mathrm{d}t}{\mathrm{d}\tau} = \sigma(\mathbf{p}, \mathbf{q})
,\end{equation}
where \(\sigma > 0\) is the step size control function. To ensure that the transformed system remains Hamiltonian, the Darboux-Sundman transformation \citep{nacozy1977intermediate} is expressed as follows:
\begin{equation}
    K = \sigma(H - H_0)
.\end{equation}
Here, \(H\) is the Hamiltonian of the original system, \(K\) is the new Hamiltonian system, and \(H_0 = H(\mathbf{p}_0,\mathbf{q}_0)\) is the initial energy of the system. For the N-body problem, the Hamiltonian system is generally separable. However, introducing the time transformation makes the system nonseparable. As a result, explicit methods become implicit, requiring the solution of nonlinear equations. Furthermore, due to the time transformation, the Jacobian matrix of the new Hamiltonian system becomes more complicated, which significantly increases the computational cost if solving nonlinear equations using Newton's method. An efficient quasi-Newton method based on Broyden's approach for solving nonlinear equations is presented in Section \ref{qnbaesdbm}.

In some cases, appropriate step size control functions \(\sigma(\mathbf{p}, \mathbf{q})\) are known a priori. For example, \citet{budd2003handbook} proposed an effective step size control function for the two-body problem:
\begin{equation}
    \sigma_1(\mathbf{p}, \mathbf{q}) = \|\mathbf{q}\|^\alpha
    \label{tb}
,\end{equation}  
where \(\alpha = 2\) or \(1.5\). Smaller steps are taken when the bodies are close to each other. In the \(N\)-body problem, close encounters can lead to a significant increase in velocity, when the total energy of the system is primarily dominated by kinetic energy and a smaller time step should be chosen in such cases. Therefore, it is natural to use the following step size control function \citep{huang1997adaptive,examplebook}:
\begin{equation}
    \sigma_2(\mathbf{q}) = \left[ \left( H_0 - U(\mathbf{q}) \right) + \nabla U(\mathbf{q})^T M^{-1} \nabla U(\mathbf{q}) \right]^{-1/2}
    \label{sigma}
,\end{equation}
where \(U(\mathbf{q})\) is the potential energy. Equation~\ref{sigma} ensures that the step size adapts naturally to the dynamic properties of the system. To prevent the step size from becoming excessively small or large, and upper and lower bounds are introduced \citep{huang1997adaptive}. The step size function used in this paper is defined as:
\begin{equation}
    {\sigma}(p, q) = \frac{\sqrt{{\sigma_2}^2 + a^2}}{\frac{1}{b} \sqrt{{\sigma_2}^2 + a^2} + 1}
    \label{sigmar}
,\end{equation}  
where \(0 < a \ll 1 \ll b\). This formulation ensures that \({\sigma}(p, q)\) is constrained within the range defined by the lower bound \(a\) and the upper bound \(b\). In this paper, \({\sigma}(p, q)\) is adopted as the time step function with \(a = 10^{-6}\) and \(b = 10^{2}\).

\subsection{An efficient quasi-Newton method based on Broyden's method}
\label{qnbaesdbm}
We let \(\mathbf{F}: \mathbf{R}^d \to \mathbf{R}^d\) be a vector valued function, defined by Eq.~\ref{Eq:14a}. Then, the Newton method for solving the nonlinear equations \(\mathbf{F}(\mathbf{x}) = 0\) can be expressed as follows:
\begin{equation}
   {{\bf{x}}_{k + 1}} = {{\bf{x}}_k} - {\bf{J}}_k^{ - 1}{\bf{F}}\left( {{{\bf{x}}_k}} \right)
,\end{equation}
where $k$ is the iteration index. In each iteration of solving the nonlinear equations, the Jacobian matrix, ${\bf{J}}_k$, and its inverse, ${\bf{J}}_k^{ - 1}$, are typically required. However, deriving the Jacobian analytically becomes highly complex due to the time transformation applied to the Hamiltonian. A common method to approximate the Jacobian is numerical differentiation. While effective, this approach can be computationally expensive, especially when evaluating the nonlinear function, ${\bf{F}}$, involves Gaussian quadrature.

To accelerate the computation, we employ the quasi-Newton method instead of standard Newton method. Specifically, we used Broyden's method, which iteratively updates the inverse of the Jacobian matrix without explicitly computing the Jacobian itself. This allowed us to significantly reduce the computational overhead while maintaining reasonable convergence properties.

Following by Broyden's method, the inverse Jacobian at iteration \( k+1 \) was updated using information from the previous iteration, \( k \), according to the following update formula:
\begin{equation}
\mathbf{J}_{k + 1}^{-1} = \mathbf{J}_k^{-1} + \frac{\left( \mathbf{s}_k - \mathbf{J}_k^{-1} \mathbf{y}_k \right) \mathbf{s}_k^\top \mathbf{J}_k^{-1}}{\mathbf{s}_k^\top \mathbf{J}_k^{-1} \mathbf{y}_k}
\label{eq:broyden_update}
,\end{equation}  
where \( \mathbf{s}_k = \mathbf{x}_k - \mathbf{x}_{k-1} \) and \( \mathbf{y}_k = \mathbf{F}(\mathbf{x}_k) - \mathbf{F}(\mathbf{x}_{k-1}) \). This formula provides an efficient method for updating \( \mathbf{J}_k^{-1} \), without requiring a direct recalculation of the Jacobian matrix. This is especially beneficial in cases where \( \mathbf{F}(\mathbf{x}) \) involves computationally intensive procedures, such as Gaussian integration. The pseudocode used in this study for the quasi-Newton method based on the Broyden method, is given in Algorithm \ref{algo:nonlinear_solver}, presented below.
\begin{algorithm}[!h]
\caption{A quasi-Newton method based on Broyden's method}\label{algo:nonlinear_solver}
\begin{algorithmic}[1]
\Require \( \mathbf{x}_0 \) (initial guess), \( \mathbf{J}_0^{-1} \) (initial inverse of the Jacobian matrix), \( \varepsilon \) (convergence criterion), \( k_{\max} \) (maximum number of iterations)
\Ensure \( \mathbf{x}_k \) (solution)
\State Initialize \( k \gets 0 \), \( e_k \gets 1, \)
\While{\( e_k > \varepsilon \) \textbf{and} \( k < k_{\max} \)},
    \State Evaluate nonlinear function: \( \mathbf{F}(\mathbf{x}_k), \)
    \State Update the solution: \( \mathbf{x}_{k+1} \gets \mathbf{x}_k - \mathbf{J}_k^{-1} \mathbf{F}(\mathbf{x}_k), \)
    \State Update the inverse Jacobian \( \mathbf{J}_{k+1}^{-1} \) using the Broyden update formula (Eq.~\ref{eq:broyden_update}),
    \State Compute the relative error: \( e_{k+1} \gets \frac{\|\mathbf{x}_{k+1} - \mathbf{x}_k\|}{\|\mathbf{x}_k\|} \),
    \State Increment the iteration counter: \( k \gets k + 1, \)
\EndWhile,
\State \Return \( \mathbf{x}_k \).
\end{algorithmic}
\end{algorithm}

Throughout the run of Algorithm \ref{algo:nonlinear_solver}, the solution from the previous time step \( \mathbf{x}_k \) is utilized as the initial guess \( \mathbf{x}_0 \) for the subsequent time step. Similarly, the inverse Jacobian, \( \mathbf{J}_0^{-1} \), is initialized by the inverse Jacobian, \( \mathbf{J}_k^{-1} \), from the prior converged time step. For the first time step, the Jacobian \( \mathbf{J}_0 \) is approximated numerically using the difference quotient method.

Nonlinear equations can also be solved by using other solvers such as the tensor method. The tensor method introduces an additional tensor term in the iteration to capture the higher-order nonlinear behavior \citep{schnabel1984tensor}. Compared to quasi-Newton methods (e.g. Broyden’s method, as used in this work), the tensor method comes with higher computational costs \citep{nocedal2006numerical} and is generally more difficult to implement due to the introduction of the additional tensor term.

However, one notable problem of the quasi-Newton method is its potentially slower convergence compared to the standard Newton iteration \citep{broyden1970methods}. This issue can often be mitigated by employing smaller step sizes, which improve the stability and accuracy of the iterations while maintaining the efficiency of the quasi-Newton approach \citep{bulin2021efficient,nocedal2006numerical}.
\section{Numerical example}\label{sec:ne}
To evaluate the performance of the proposed integration methods,
we present three numerical examples. Before delving into the examples, we summarize the key characteristics of the integrators involved.

The original SQQ uses equidistant interpolation nodes and solves nonlinear equations with the Newton method. Five extensions of SQQ are evaluated in this Section: SQQ-P, SQQ-PN, SQQ-PQ, SQQ-PTN, and SQQ-PTQ. The suffix "T" denotes the use of the time transformation, as defined in Eq.~\ref{sigmar}, "P" refers to the projection method for efficient interpolation, while "Q" and "N" indicate the use of quasi-Newton and Newton methods for solving nonlinear equations, respectively. The characteristics of these methods are summarized in Table~\ref{tab:algorithms}. For consistency, the convergence threshold, \( \varepsilon \), for solving nonlinear equations is set to $1 \times 10^{-12}$ for all examples. All numerical computations were performed in MATLAB on a platform equipped with a 13th Generation Intel(R) Core(TM) i7-13700 processor and 32~GB of RAM.
\begin{table*}
    \centering
    \fontsize{8}{9.6}\selectfont 
    \renewcommand{\arraystretch}{1.2}
    \caption{Summary of characteristics for different integrators}
    \label{tab:algorithms}
    \setlength{\tabcolsep}{10pt}
    \begin{tabular}{lcccc}
        \toprule
        Integrator name & Interpolation method & Projection method & Time transformation & Nonlinear solver \\
        \midrule
        SQQ        & Equidistant & No  & No  & Newton method \\
        SQQ-P      & Chebyshev   & Yes & No  & Newton method \\
        SQQ-PN     & Chebyshev   & Yes & No  & Newton method \\
        SQQ-PQ     & Chebyshev   & Yes & No  & Quasi-Newton method \\
        SQQ-PTN    & Chebyshev   & Yes & Yes & Newton method \\
        SQQ-PTQ    & Chebyshev   & Yes & Yes & Quasi-Newton method \\
        \bottomrule
    \end{tabular}
\end{table*}

\subsection{The Kepler problem} \label{twobody}
The Kepler problem with the Hamiltonian is given by
\begin{equation}
H = \frac{1}{2} (p_1^2 + p_2^2) - \frac{1}{r}, \quad r = \sqrt{q_1^2 + q_2^2}
,\end{equation}  
where the initial conditions are
\begin{equation}  
\mathbf{q} = \begin{bmatrix} 1-e \\ 0 \end{bmatrix}, \quad  
\mathbf{p} = \begin{bmatrix} 0 \\ \sqrt{\frac{1+e}{1-e}} \end{bmatrix},
\end{equation} 
and \( e \) represents the eccentricity of the orbit. The period of this system is \( 2\pi \).

To evaluate the efficiency of the projection method, $e = 0.5$ is selected as the test parameter. Simulations of 500 orbital periods are performed using both the SQQ and SQQ-P.

\begin{figure}[!ht]
\centering
\includegraphics[width=\linewidth]{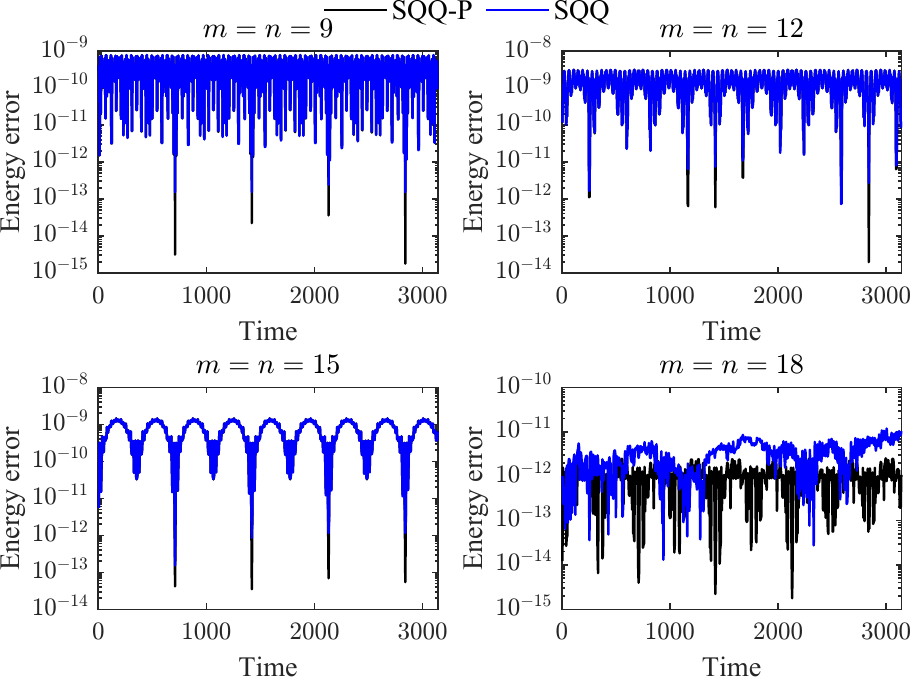}
\caption{Evolution of the energy error for SQQ-P and SQQ. Here, $e = 0.5$, and $m$ and $n$ denote the number of interpolation points for the generalized displacement and generalized momentum, respectively.}
\label{tb_error1}
\end{figure}

Figure~\ref{tb_error1} illustrates the evolutions of the energy error for SQQ and SQQ-P. When the number of interpolation points is relatively small (i.e., $m < 15$ and $n < 15$), the Runge phenomenon does not occur. In this example, the evolutions of the energy error for SQQ and SQQ-P are nearly identical, exhibiting oscillatory behavior within a bounded range, which is typical for symplectic methods. This demonstrates that the projection method employed in this study provides an accurate computation of the interpolation function. As the number of interpolation points increases (i.e., $m = 18$ and $n = 18$), SQQ exhibits a growth in energy error over time. This behavior arises from the excessive number of interpolation points, which triggers the Runge phenomenon. As a result, the interpolated action, $S$, may become inaccurate. In contrast, SQQ-P, which uses Chebyshev interpolation points, effectively suppresses the Runge phenomenon. As a result, the energy error for SQQ-P continues to oscillate within a bounded range, as shown in Fig.~\ref{tb_error1}. This highlights the robustness of Chebyshev interpolation in maintaining the stability of the computed action.

Table \ref{tab:cputInTb} shows the time step sizes and CPU times for the various simulations. As the number of interpolation points increases, the CPU time ratio decreases (the CPU time ratio is defined as the ratio of SQQ-P's CPU time to that of SQQ). This happens because the increase in the number of interpolation points introduces additional Gaussian nodes, significantly raising the computational cost of evaluating the interpolation function during each time step. In contrast, SQQ-P employs the projection method proposed in this study,  thereby enhancing the efficiency.
\begin{table}[h]
\centering
\caption{CPU time comparison between SQQ and SQQ-P methods.}
\label{tab:cputInTb}
\fontsize{9}{11}\selectfont
\begin{tabular}{lcccc}
\toprule
$m, n$ & 9 & 12 & 15 & 18 \\
\midrule
$\Delta t$ & 0.4 & 0.6 & 0.8 & 1.0 \\
SQQ CPU time (s) & 54.77 & 55.07 & 67.78 & 113.01 \\
SQQ-P CPU time (s) & 29.82 & 27.77 & 35.06 & 51.51 \\
CPU time ratio (\%) & 55.46 & 50.43 & 51.72 & 45.58 \\
\bottomrule
\end{tabular}
\tablefoot{
CPU time ratio is defined as the ratio of SQQ-P's CPU time to that of SQQ expressed as a percentage.
}
\end{table}

\begin{figure}[!h]
\centering
\includegraphics[width=\linewidth]{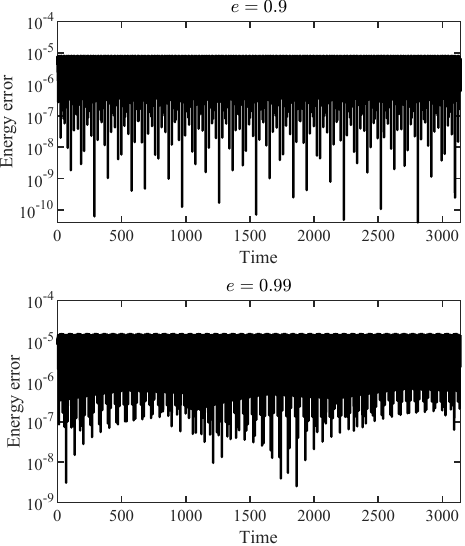}
\caption{Evolution of the energy error for SQQ-PTQ for the Kepler problem with high eccentricities.}
\label{SQQPTQERROR_twobody_500P}
\end{figure}

\begin{figure}[!ht]
\centering
\includegraphics[width=\linewidth]{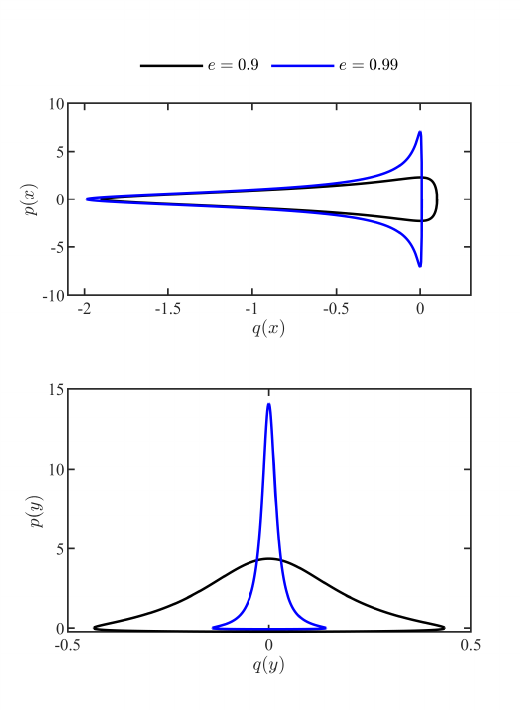}
    \caption{Phase space trajectories over 500 orbital periods for the Kepler problem with different eccentricities.}
\label{tb_phase}
\end{figure}

The Kepler problem becomes challenging for traditional fixed-step algorithms as the eccentricity approaches to 1. We perform two numerical experiments for eccentricities of $e = 0.9$ and $0.99$, with parameters $\Delta \tau = 0.01$, $m = n = 3$, and a total integration time of 500 orbital periods. Fig.~\ref{SQQPTQERROR_twobody_500P} illustrates the 
evolution of the energy error for SQQ-PTQ with these large eccentricities. Fig.~\ref{tb_phase} shows the corresponding phase space trajectories of \( q \) and \( p \). It can be seen that, even over long-term integration, the phase-space structure remains nearly intact, which indicates that SQQ-PTQ effectively preserves the Hamiltonian structure.

Table \ref{tab:performanceforhe2body} presents the CPU times of SQQ-PTN and SQQ-PTQ for solving the Kepler problem with high eccentricities. SQQ-PTQ exhibits an improvement in efficiency over SQQ-PTN, with the latter requiring about 2 to 3 times of the CPU time at $e = 0.9$ and $e = 0.99$, respectively.

\begin{table}[h]
\centering
\fontsize{9}{10.8}\selectfont
\caption{Performance comparison of SQQ-PTN and SQQ-PTQ for the Kepler problem at high eccentricities.}
\label{tab:performanceforhe2body}
\begin{tabular}{lcc}
\toprule
Integrator   & SQQ-PTN & SQQ-PTQ  \\
\midrule
CPU time (s) for $e = 0.9$  & 231.63 & 98.99 \\
CPU time (s) for $e = 0.99$ & 898.76 & 282.63 \\
\bottomrule
\end{tabular}
\end{table}

\subsection{Three-body problem} \label{subsec_three-body}
\begin{figure}[!h]
\centering
\includegraphics[width=\linewidth]{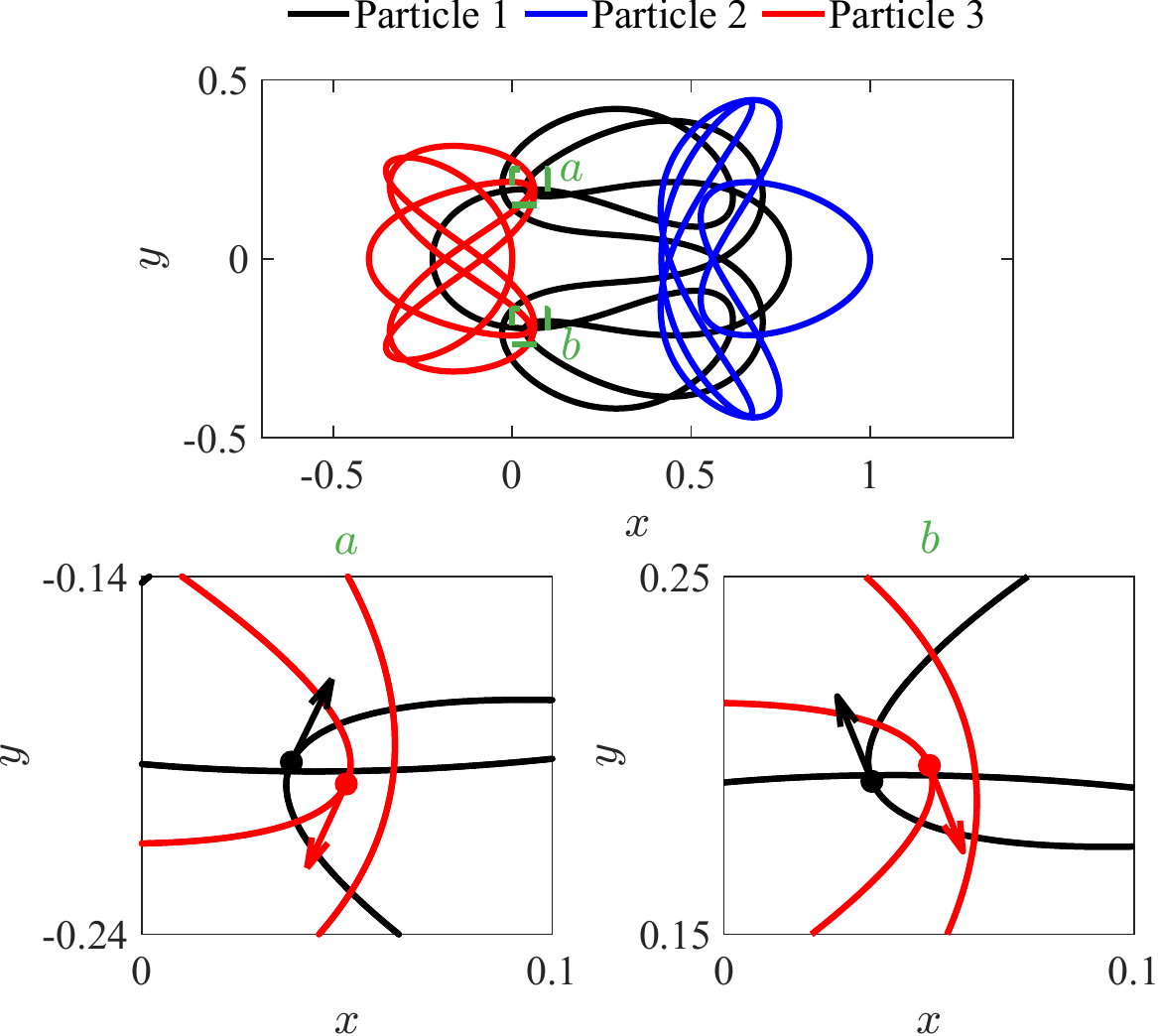}
\caption{Trajectories of particles in the three-body problem, with dashed boxes $a$ and $b$ indicating the locations of two close encounters. Box $a$ corresponds to a close encounter occurring at approximately $0.4T$, and box $b$ corresponds to the second close encounter occurring at approximately $0.6T$, where $T$ denotes the orbital period. The two subplots below provide zoomed-in views, showing the positions and velocity directions of particles 1 and 3. For both encounters, the closest distance between particles 1 and 3 is approximately 0.014.}
\label{fig:threebodytra}
\end{figure}

\begin{figure}[!h]
\centering
\includegraphics[width=\linewidth]{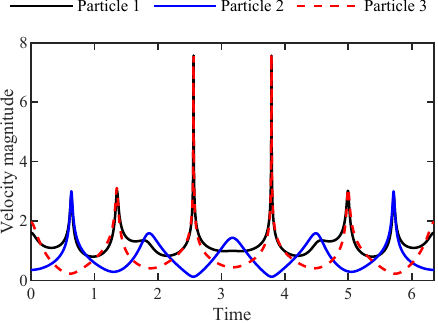}
\caption{Evolution of the velocity magnitudes of the three particles over one period in the three-body problem.}
\label{fig:threebodyvel_magn}
\end{figure}

\begin{figure}[!h]
\centering
\includegraphics[width=\linewidth]{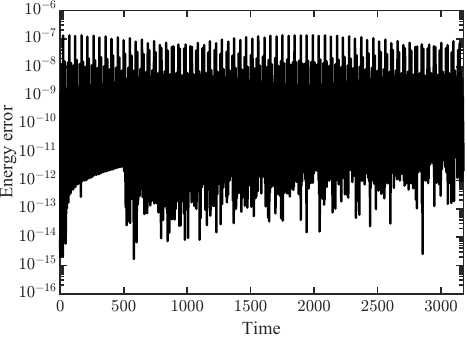}
\caption{Evolution of the energy error for the integration of the three-body problem by SQQ-PTQ.}
\label{fig:threebodyerror}
\end{figure}

\begin{figure}[!h]
\centering
\includegraphics[width=\linewidth]{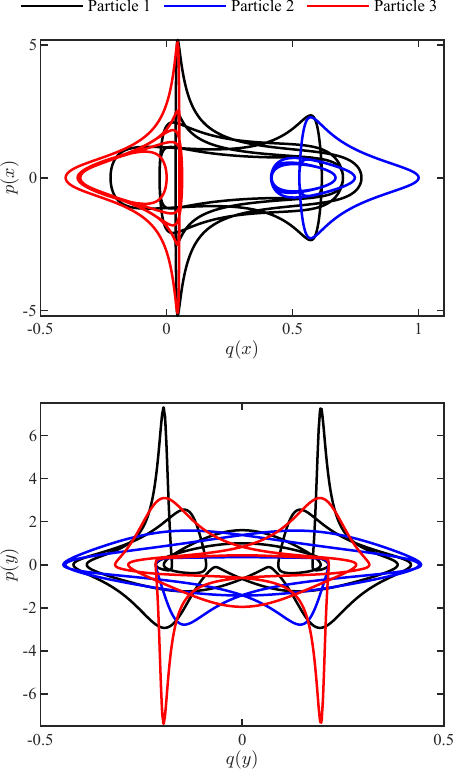}
\caption{Phase space trajectories of particles in the three-body problem over 500 orbital periods.}
\label{fig:threebodyPhase}
\end{figure}

For problems involving close encounters, any unreasonable time step will lead to increased errors or calculation failure, so an appropriate time step is necessary to accurately capture such events. The The three-body problem with possible close encounters is a good option for assessing the performance of the adaptive time step strategy of integrators. Recently, \citet{liao2022three} employed machine learning techniques to discover some novel periodic solutions for the three-body problem. In this section, we use one specific periodic configuration from their study to test SQQ-PTQ's ability to automatically select the time step. The orbital and physical parameters of the three bodies, as given in \citet{liao2022three}, are as follows:
\begin{equation}
\begin{array}{lll}
\mathbf{q}_1 = \left[ -0.2227,\ 0 \right]^{\mathrm{T}}, & 
\dot{\mathbf{q}}_1 = \left[ 0,\ 1.7813 \right]^{\mathrm{T}}, & 
m_1 = 0.9, \\
\mathbf{q}_2 = \left[ 1,\ 0 \right]^{\mathrm{T}}, & 
\dot{\mathbf{q}}_2 = \left[ 0,\ 0.4150 \right]^{\mathrm{T}}, & 
m_2 = 0.85, \\
\mathbf{q}_3 = \left[ 0,\ 0 \right]^{\mathrm{T}}, & 
\dot{\mathbf{q}}_3 = \left[ 0,\ -1.9559 \right]^{\mathrm{T}}, & 
m_3 = 1.
\end{array}
\label{Eq: 39}
\end{equation}
The period of the system is ${T} = 6.3509$, and the gravitational constant $G$ is normalized to 1. We simulated 500 orbital periods of this system using the SQQ-PTQ. Figure~\ref{fig:threebodytra} shows the trajectories of three particles, where particles 1 and 3 experience two close encounters within one period. These encounters occur at approximately $0.4T$ (dashed box $a$), and at $0.6T$ (dashed box $b$). For both encounters, the closest distance between particles 1 and 3 is approximately 0.014. Figure~\ref{fig:threebodyvel_magn} shows the magnitude of the velocities of the three particles (during the first period). The velocities of particles 1 and 3 exhibit two peaks within one period, due to their close encounters, leading to significant velocity changes for both particles. These two close encounters cause the maximum velocities of particles 1 and 3 to be 10 and 34 times of their respective minimum velocities. The evolution of the energy error exhibits stable oscillations as shown in Fig.~\ref{fig:threebodyerror}. Furthermore, the phase space trajectories shown in Fig.~\ref{fig:threebodyPhase} maintain a well-preserved structure throughout long-term integration, indicating that SQQ-PTQ effectively preserves the symplectic property.

We also use the MATLAB built-in general-purpose integrator ODE45 (with relative error tolerance set to $1 \times 10^{-10}$) and SQQ-PTN to simulate this system for 500 orbital periods. Table \ref{tab:performancefor3body} lists the computational time and maximum energy errors for the three integrators. The only difference between SQQ-PTQ and SQQ-PTN lies in the method used to solve the nonlinear equations, which explains why their maximum energy errors are nearly identical. The CPU time of SQQ-PTQ is less than one-third of that of SQQ-PTN. This indicates that the quasi-Newton method accelerates the computation without introducing any additional errors.

\begin{table}[h]
\centering
\fontsize{9}{10.8}\selectfont
\caption{Comparison of the performances between SQQ-PTN, SQQ-PTQ, and ODE45 (with relative error tolerance set to $1 \times 10^{-10}$) in the three-body problem.}
\begin{tabular}{lccc}
\toprule
Integrator       & SQQ-PTN & SQQ-PTQ & ODE45  \\
\midrule
$m, n$            & 3       & 3       & -      \\
$\Delta\tau$     & 0.01       & 0.01       & -      \\
CPU time (s)     & 137.04  & 43.62   & 87.17  \\
Maximum energy error & 1.32E-7 & 1.32E-7 & 3.83E-7 \\
\bottomrule
\end{tabular}
\label{tab:performancefor3body}
\end{table}


\subsection{Outer Solar System}
As the final example, we consider the outer Solar System problem, which includes the Sun, Jupiter, Saturn, Uranus, Neptune, and Pluto. The outer Solar System has multiple degrees of freedom. Its clear structure and stable orbits over long timescales make it suitable for evaluating the integrator’s robustness, stability, and computational efficiency. For a detailed description of the masses and orbital parameters of the planets, we refer to \citet{examplebook}. 

\begin{figure}[!ht]
\centering
\includegraphics[width=\linewidth]{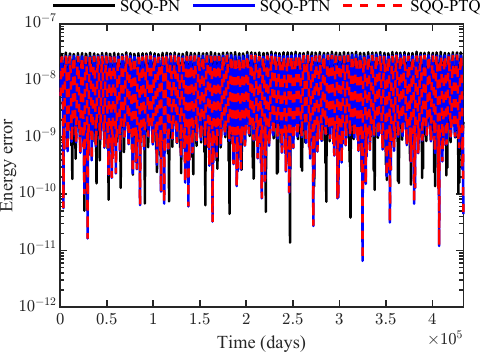}
\caption{Evolutions of the energy error for the integration of the outer Solar System by SQQ-PN, SQQ-PTN, and SQQ-PTQ, with a total integration time of 100 Jupiter orbital periods.}
\label{os-error}
\end{figure}

We perform simulations over 100 Jupiter orbits using three integrators: SQQ-PN, SQQ-PTN, and SQQ-PTQ. Figure~\ref{os-error} illustrates the evolution of energy error over time. It is worth mentioning that the evolutions of the energy error for SQQ-PTQ and SQQ-PTN are nearly indistinguishable. This is consistent with the phenomenon of the three-body problem mentioned in Section \ref{subsec_three-body}, because the only difference between SQQ-PTQ and SQQ-PTN lies in the method used to solve the nonlinear equations. Table \ref{tab:performanceos} presents a summary of the CPU time and maximum energy errors for the three integrators. It can be seen that SQQ-PTN has the longest computational time, which is roughly six times that of SQQ-PTQ. This is primarily due to the time transformation, which increases the computational cost of the Jacobian matrix in each iteration. SQQ-PN, which does not use the time transformation, exhibits moderate computational efficiency. In contrast, SQQ-PTQ achieves the highest efficiency by using a quasi-Newton method based on Broyden's method, avoiding the direct computation of the Jacobian matrix. Even compared to SQQ-PN, SQQ-PTQ requires only 50\% of the computational time.

\begin{figure}[!ht]
\centering
\includegraphics[width=\linewidth]{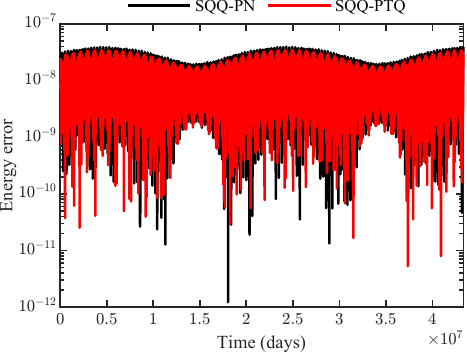}
\caption{Evolutions of the energy error for the integration of the outer Solar System by SQQ-PN and SQQ-PTQ, with a total integration time of 10000 Jupiter orbital periods.}
\label{os-erro-10000r}
\end{figure}

To further evaluate the long-term integration stability of SQQ-PTQ, we performed simulations over 10,000 Jupiter orbital periods using both SQQ-PTQ and SQQ-PN. Figure~\ref{os-erro-10000r} shows that the evolutions of the energy error for both SQQ-PN and SQQ-PTQ oscillate within a bounded range. Table \ref{tab:longperformanceos} shows the CPU time and maximum energy error of the integration with different integrators. The efficiency of SQQ-PTQ is higher than that of SQQ-PN, even though SQQ-PTQ includes the time transformation. This demonstrates the effectiveness and robustness of combining the time transformation with the quasi-Newton method.

\begin{table}[h]
\centering
\fontsize{8.5}{10.2}\selectfont
\setlength{\tabcolsep}{4pt} 
\caption{Maximum energy error and CPU times for the integration of the outer Solar System by SQQ-PN, SQQ-PTN, SQQ-PTQ, and ODE45 (relative error tolerance $1 \times 10^{-8}$), over 100 Jupiter orbital periods.}
\label{tab:performanceos}
\begin{tabular}{lcccc}
\toprule
Integrator & SQQ-PTN & SQQ-PN & SQQ-PTQ & ODE45 \\
\midrule
$m, n$ & 5 & 5 & 5 & -- \\
$\Delta \tau$ (days) & 250 & 250 & 250 & -- \\
CPU times (s) & \num{68.91} & \num{25.92} & \num{11.01} & \num{50.33} \\
Max energy error & \num{2.85e-8} & \num{3.21e-8} & \num{2.85e-8} & \num{4.71e-7} \\
\bottomrule
\end{tabular}
\end{table}

\begin{table}[ht]
\centering
\fontsize{9}{10.8}\selectfont
\caption{Maximum energy error and CPU times for the integration of the outer Solar System by SQQ-PN, SQQ-PTQ and ODE45 (with relative error tolerance set to $1 \times 10^{-8}$), with a total integration time of 10000 Jupiter orbital periods.}
\label{tab:longperformanceos}
\begin{tabular}{lccc}
\toprule
Integrator  & SQQ-PN & SQQ-PTQ & ODE45 \\
\midrule
$m, n$ & 5 & 5 & -- \\
$\Delta \tau$ (days)  & 250 & 250 & -- \\
CPU times ($\times 1000$ s)  & 2.26 & 1.31 & 3.56 \\
Maximum energy error  & 3.97E-08 & 3.51E-08 & 3.13E-04 \\
\bottomrule
\end{tabular}
\end{table}

\section{Conclusion}
This paper introduces an adaptive symplectic integrator SQQ-PTQ, which builds on the fixed-step integrator SQQ. By combining the time transformation with a quasi-Newton method, SQQ-PTQ improves both computational efficiency and applicability.

A key feature of SQQ-PTQ is its use of time transformation to achieve adaptive step sizing. This allows SQQ-PTQ to adjust the step size in gravitational multi-body problems on different timescales. To overcome the challenges of computing Jacobian matrices, SQQ-PTQ employs a quasi-Newton method based on Broyden's method. This approach significantly accelerates the solution of nonlinear equations, enhancing the overall computational performance.

The effectiveness and robustness of SQQ-PTQ are demonstrated through three numerical experiments. In particular, SQQ-PTQ demonstrates adaptability when handling close-encounter problems. Furthermore, over long-term integrations, SQQ-PTQ maintains energy conservation, further confirming its characteristics as a symplectic algorithm.

\begin{acknowledgements}
This work was supported by the National Natural Science Foundation of China (No.~12311530055 and 12472048), the Shenzhen Science and Technology Program (Grant No.~ZDSYS20210623091808026), and by High-performance Computing Public Platform (Shenzhen Campus) of SUN YAT-SEN UNIVERSITY.
\end{acknowledgements}

\end{document}